\documentclass[12pt,twoside,a4paper]{amsart}
\usepackage{amsmath,amssymb,amsthm}
\usepackage{stackrel}
\usepackage{setspace}
\usepackage{lscape}
\usepackage{multirow}
\usepackage{rotating}
\usepackage{graphicx}
\usepackage{hyperref}
\usepackage{hyperref}
\hypersetup{
    colorlinks=true,
    linkcolor=blue,
    filecolor=magenta,      
    urlcolor=cyan,
    citecolor=blue}
\calclayout    

\newtheorem{Thm}{Theorem}[section]
\newtheorem{Lem}[Thm]{Lemma}
\newtheorem{Pro}[Thm]{Proposition}
\newtheorem{Cor}[Thm]{Corollary}

\newtheorem*{acknowledgement*}{Acknowledgement}

\theoremstyle{definition}

\theoremstyle{remark}
\newtheorem{Rem}[Thm]{Remark}
\newtheorem{Def}[Thm]{Definition}

\newcommand{\R}{\mathbb{R}}

\renewcommand{\phi}{\varphi}

\newcommand{\Vol}{\operatorname{Vol}}

\newcommand{\tr}{\operatorname{tr}}

\renewcommand{\cos}{\operatorname{cos}}

\begin{document}
\title[On Harmonic and Asymptotically Harmonic Finsler Manifolds]{On Harmonic and Asymptotically Harmonic Finsler Manifolds}
\author {Hemangi  Shah and Ebtsam  H. Taha}
\address{Harish-Chandra Research Institute, Chhatnag Road, Jhunsi, Allahabad 211019, India}
\email{hemangimshah@hri.res.in}
\vspace{-1 cm}
\address{Department of Mathematics, Faculty of Science, Cairo University, Giza 12613, Egypt}
\vspace{-1 cm}
\address{Harish-Chandra Research Institute, Chhatnag Road, Jhunsi, Allahabad 211019, India}
\email{ebtsam.taha@sci.cu.edu.eg}
\thanks{ This paper has been presented in the  international workshop which held in honor of Misha Gromov \lq \lq 2020 Virtual Workshop on Ricci and Scalar Curvature"}
\maketitle 
\begin{abstract}In the present paper we introduce and investigate various types of harmonic Finsler manifolds and find out the interrelation between them. We give some characterizations of such spaces in terms of the mean curvature of geodesic spheres and the Laplacian of the distance function induced by the Finsler structure. We investigate some properties of the Finsler mean curvature of geodesic spheres of different radii. In addition, we prove that certain harmonic Finsler manifolds are of Einstein type and provide a technique to construct harmonic Finsler manifolds of Randers type. Moreover, we give some examples of non-Riemmanian Finsler harmonic manifolds of constant flag curvature and constant $S$-curvature.  
\end{abstract}
\maketitle 
\bigskip


\maketitle 

\section{Introduction}
The notion of harmonic manifolds was introduced by H. S. Ruse \cite{Ruse}  in  the first half of the $20^{th}$ century. As rich interplay of analysis and geometry continued to grow so did the fertility of these sub-disciplines of differential geometry. In fact, the study of  harmonic and asymptotically harmonic Riemannian manifolds is an active area of research. A complete Riemannian manifold M is said to be \textit{harmonic} if all geodesic
spheres in $M$ of sufficiently small radii are of constant mean curvature.  Several results and characterizations of harmonic manifolds appeared in  \cite{Cay,  harmonicbook16, Besse, RamRanj, Shah03, {Shahth}}. For example, harmonic Riemannian manifolds can be characterized by the mean value property of harmonic functions which was proved by Willmore \cite{Willmore}. It is known that harmonic Riemmanian manifolds have constant Ricci curvature so that they are Enistein manifolds \cite{Besse}. A complete classification of compact harmonic manifolds had been achieved by Szabo in \cite{Szabo}, in which he proved the Lichnerowicz conjecture \lq \lq \textit{every compact simply connected harmonic manifold is either flat or rank one symmetric space}". Recently, in \cite{new-har}, the authors discovered a new class of harmonic Hadamard manifolds and established a characterization of  harmonic Hadamard manifolds of hypergeometric type with respect to the volume density.\\
 
 Finsler geometry  is a further generalization of Riemannian geometry and is much wider in scope and richer in content. For instance, the model spaces (space forms) in Riemannian geometry \cite{ Pet16} are well understood
and classified; however, in Finsler geometry the problem is far from being completely solved. Some partial results, for example \cite{BCS, {Shen1996}, Shenfish}, in certain special Finsler spaces indeed exist.
In fact, there are infinitely many Finsler
model spaces, which are not isometric to each other. This difficulty persists in Finsler geometry even with the special cases  like that of constant flag curvature, due to the abundance of geometric objects associated with the Finsler structure like Cartan torsion, Berwald curvature, Dougluas tensor, S-curvature, T-curvature,... etc,  which all vanish identically in the Riemannian case, cf.~\cite{BCS,  OhtaSzero, curdisvol, Shlec,  ShenGeometricmeaning,  Shenbook16, ShiBanktesh}. 
Indeed, working in the Finsler context needs different techniques that do not exist in the Reimannian case. This makes the  study of Finsler problems challenging.\\

  One of the central focus of study is the Riemann-Finsler geometry, that is, the area where geometers are interested in generalizing Riemannian results to the Finsler context. In this direction, we generalize harmonic manifolds to the Finsler case. Such generalizations have not been  studied in the literature before and are inspired by \cite{Besse}. First we introduce several types harmonic manifolds in the Finsler context, viz. locally, globally, infinitesimal, asymptotic harmonic Finsler manifolds. 
 Our results on harmonic Finsler manifolds reduce to the Riemmanian ones when the Finsler metric is Riemmanian. To the best of our knowledge, the only papers dealing with harmonic Finsler manifolds are \cite{Kim} and \cite{Kim2}. However, our treatment and results both are completely different. \\
 
 We study the relations between these notions and, in particular, prove that certain harmonic Finsler manifolds are of Einstein type. Different characterizations of such spaces are established in terms of Shen's Finsler Laplacian, Finsler mean curvature and isoparametric Finsler distance.  To enrich understanding, we provide various examples of non-Riemmanian Finsler harmonic manifolds such as  \textit{Minkoskian} metrics, \textit{Funk} metrics, \textit{Shen's fish tank} metric and a family of non-Riemannian Finsler metrics on odd-dimensional spheres constructed by Bao and Shen. The first two metrics are projectively flat whereas the last two metrics are not. Additionally, we  give a technique to construct  harmonic Finsler manifolds of Randers type in Theorem~\ref{harmonicRanderthm}. \\

 In what follows, we give the structure of the present paper. Section 2 is devoted to some preliminaries needed for better exposition of our results. Thereafter, in next four sections we present our main results. 
 In particular, in section 3 we study some properties of normal and mean curvatures of geodesic spheres in Finsler manifolds.
Section 4  deals with the formulation of various types of harmonic Finsler manifolds and their inter-relationships as well as some examples of such spaces. Moreover, we do study compact and noncompact harmonic manifolds. In section 5, we investigate harmonic Finsler manifolds of Randers type.

 \section{Preliminaries}
We will use the following notations.  $M$ denotes an $n$-dimensional, $n>1$, smooth connected orientable manifold, 
$(TM, \pi, M)$, or simply $TM$, its tangent bundle and $TM_{0}:=TM\setminus \{0\}$ the tangent bundle with the null section removed. The tangent vector space at each $x \in M$ without the zero vector is denoted by  $T_x M_{0}$.  The local coordinates $(x^{i})$ on $M$ induce  local
coordinates $(x^{i},y^{i})$ on $TM$.  The pullback bundle of $TM$ is  denoted by $\pi^{-1}(TM)$. Moreover, 
$\partial_i$ denotes partial differentiation with respect to $x^i$ and  $\dot{\partial}_i$ denotes partial differentiation
    with respect to  $y^i$ (basis vector fields of the vertical bundle). 
 
\begin{Def} \label{Finslerdef}\cite{BCS} A \emph{smooth Finsler structure} on a manifold $M$  is a mapping $F:TM\rightarrow [0, \infty )$ with the following properties:
\begin{itemize}
    \item[(a)] $F$ is $C^\infty$ on the slit tangent bundle  $TM_0$.

    \item[(b)] $F$ is positively homogeneous of degree one in $y$: $F(x,\lambda y) = \lambda F(x,y)$
     for all $y \in T_{x}M$ and $\lambda > 0$.

    \item[(c)] The Hessian matrix $(g_{ij}(x,y))_{1 \leq i,j \leq n}$
is positive definite at each point $y$ of $TM_0$, where $\displaystyle{g_{ij}(x,y):=\frac{1}{2} \dot{\partial}_i\dot{\partial}_j F^2(x,y)}$.
\end{itemize}
The pair $(M,F)$ is called  a \emph{Finsler manifold} and the symmetric bilinear form $g_{ij}(x,y)$ is called the
\emph{Finsler metric} of the Finsler structure $F$.
\end{Def}
\begin{Rem}
$(i)$ A Finsler metric is Riemannian when  $g_{ij}(x,y)$ are functions in $x$ only and it is locally Minkoskian when $g_{ij}(x,y)$ are functions in $y$ only in some coordinate system.

$(ii)$ A Finsler metric can be characterized in any tangent space $T_{x}M$ by its unit vectors, which form a smooth strictly convex hypersurface $I_x M$ called \emph{indicatrix} at the point $x \in M$.  
 When a Finsler metric is Riemannian, this hypersurface at each point of $M$ is a Euclidean unit sphere \cite[\S 2.2.1]{Shenbook16}. The indicatrix of $F$ is $IM:=\displaystyle\bigcup_{x\in M}I_{x}M$.
\end{Rem}
If we relax condition (c) of Definition \ref{Finslerdef} to be $(g_{ij}(x,y))_{1 \leq i,j \leq n}$ is a nondegenerate  matrix  then we  deal with\textit{ pseudo} or \textit{nondegenerate  Finsler structure}. 
\begin{Def}\cite{BCS}
A Finsler manifold is said to be \emph{reversible} if  $F(x,-y)=F(x,y), \, \forall y \in T_{x}M$.
\end{Def}
\begin{Def}
  A Finsler manifold $(M,F)$ is of Randers type if $F:=\alpha + \beta,$ where $\alpha:= \sqrt{\alpha_{ij}(x) \,y^{i}\, y^{j}}$ \footnote{Hereafter, the Einstein summation convention is in place.} is a Riemannian metric and $\beta = b_{i}(x) \,dx^{i}$ is a $1$-form on $M$ with $|| \beta ||_{\alpha} <1$. The map $F$ is then said to be a Randers metric.
\end{Def} 
 \begin{Def}  
 The coefficients of the geodesic spray of a Finsler structure  $F$ are given by \[G^i (x,v):= \frac{1}{4} g^{ij} \{y^{l}\, \partial_{l} \dot{\partial_{j}} (F^{2}) -  \partial_{j}(F^{2})\}. \] 
\end{Def}   
  Consequently,
    $ N^i_j:=\dot{\partial}_jG^i$ is a nonlinear connection; the Barthel connection
    associated with $F$, and
$\delta_i:=\partial_i-N^r_i \,\dot{\partial}_r$ are the basis
     vector fields of the horizontal bundle induced by the Barthel connection.
       
  Another special Finsler space which includes Riemannian and locally Minkoskian manifolds is the Berwald manifold. More precisely,
\begin{Def} A Finsler manifold $(M,F)$ is said to be Berwaldian if the Berwald tensor $G^{h}_{ijk}:=\dot{\partial}_i\dot{\partial}_j\dot{\partial}_kG^h$ vanishes.
\end{Def}

\textbf{Generalized metric space.} \cite{Shlec, Tam08}
 The distance $d_{F}$ induced from the Finsler structure $F$ can be defined naturally in $M$ as follows $$d_{F}(p,q):= inf \{\int^{1}_{0} F(\dot{\eta}(t))\, dt\,\ | \, \eta: [0,1]  \rightarrow M,\,  C^{1}\text{ curve joining }p \text{ to } q \}.$$
 \begin{Rem}
$(i)$ It should be noted that the \emph{Finsler distance} is nonsymmetric, that is, $d_{F}(p,q)\neq d_{F}(q,p)$. The pair $(M, d_{F})$ is sometimes called a\emph{ generalized metric space}. \\
$(ii)$  It is known that $d_{F}$ is symmetric if and only if the Finsler structure is reversible. In other words, the distance depends on the direction of curve. Therefore, the reverse of a general Finsler geodesic  
can not be a geodesic. The non-reversibility property is also reflected in the notion of Cauchy sequence and completeness \cite[\S 6.2]{BCS}.\\
$(iii)$ Thus, being different from the Riemannian case,  a positively (or forward) complete Finsler manifold $(M , F)$ is not necessarily negatively (or backward) complete. For example, a non-Riemannian Randers metric is positively complete solely. The classical Hopf-Rinow theorem splits into forward and backward versions \cite[\S 6.6]{BCS}.  A Finsler metric is called a \emph{complete} if it is both forward and backward complete.
\end{Rem} 

Another main difference between Finsler and Riemannian geometries is that in a general Finsler manifold, the exponential map is only $C^{1}$ at the origin of $T_{x} M$ (however, it is $C^{\infty}$ on $T_{x} M_{0}$).  It was proved by \textit{Akbar-Zadeh} that the exponential map is $C^{2}$ at the origin if and only if the Finsler manifold is Berwaldian \cite[\S 5.3]{BCS}. More details about the exponential map are in~\cite{Shlec, Shenbook16}. \\
 
\textbf{Volume measures in Finsler manifolds.}
\begin{Def}\cite[\S 2.1]{Shlec}
A \emph{Finsler $\mu$-space} is a Finsler manifold $(M,F)$ equipped with a volume measure $d\mu$ (nondegenerate $n$-volume form) on M. 

A volume measure $d\mu$ can be written in the local coordinates $(x^1, ..., x^n)$ as follows 
\begin{equation}\label{arbitrarysigma}
d\mu = \sigma_{\mu}(x) \,dx^{1} \wedge ... \wedge dx^{n}=\sigma_{\mu}(x) \,dx,
\end{equation} where $\sigma_{\mu}(x)$ is a positive function on $M$ satisfying certain properties.
\end{Def}
There are several, non-equivalent definitions of volume forms used in Finsler geometry.  The most well known are \emph{Busemann-Hausdorff and Holmes-Thompson volume forms} \cite{Shenbook16}. 
\begin{Def}
The \textit{Busemann-Hausdorff volume form} is defined  at a point $x \in M$, in a  local coordinate system $(x^i)$,  as follows 
\begin{equation}
dV_{BH}:= \frac{\Vol(\mathbb{B}^n(1))}{\Vol(B_{F}^n(1))}\, \,dx^{1} \wedge ... \wedge dx^{n},
\end{equation}
where $\Vol (\mathbb{B}^n(1))$ denotes the Euclidean volume of a unit Euclidean ball \cite{vlolandsberg}:
$$\Vol(\mathbb{B}^n(1))=\frac{1}{n} \Vol(\mathbb{S}^{n-1})= \frac{1}{n} \Vol(\mathbb{S}^{n-2}) \int_{0}^{\pi} \sin^{n-2}(t)\,dt, $$
 \[\Vol(B_{F}^n(1)):= \Vol \left(\{ (y^i) \in \mathbb{R}^{n}\,| \,F(x,y^i \partial_{i})<1\}\right).\]
\end{Def}
\begin{Rem} For reversible Finsler functions,  Busemann proved that the Busemann-Hausdorff volume form is the Hausdorff measure of the metric space induced by the Finsler structure \cite{Shlec}.  
\end{Rem}
\begin{Def}The {Holmes-Thompson volume form}  is defined by $dV_{HT}: = \sigma_{HT}(x) \,dx,$ where
$$\sigma_{HT}(x):=\frac{\Vol(B_{F}^n(1), g)} {\Vol(\mathbb{B}^n(1))}= \frac{1}{\Vol(\mathbb{B}^n(1))} \int_{B_{F}^n(1)} det(g_{ij}(x,y))\, dy.$$
\end{Def}

\begin{Def}\cite{volform} The maximum and minimum  volume forms for a Finsler manifold $(M, F )$ are given respectively by
$$dV_{max} = \sigma_{max}(x) \,dx^{1} \wedge ... \wedge dx^{n},\,\,\,\,dV_{min} = \sigma_{min}(x)\, dx^{1} \wedge ... \wedge dx^{n},$$ where 
\[ \sigma_{max}(x):= \max _{y \in T_{x}M_{0}} \sqrt{det(g_{ij}(x,y))},\,\,\,\,\sigma_{min}(x):= \min_{y \in T_{x}M_{0}} \sqrt{det(g_{ij}(x,y))}.\]
\end{Def}
They may be called extreme volume forms. B. Wu \cite{volform} has used them to generalize Calabi-Yau’s linear volume growth theorem.\\
 
 One of the most important geometric objects associated with the volume measure is the S-curvature. The S-curvature was introduced by Z. Shen to study volume comparison theorem in Riemann-Finsler geometry. It is connected to the flag curvature; for more details cf. \cite{Shenbook16}. 
\begin{Def}\cite{Shlec}
 The \emph{distortion} $\tau_{\mu}$ of $(M,F, \mu)$ is defined by 
 \begin{equation}\label{distortion def}
 \tau_{\mu}(x,v) := \log \left(\frac{\sqrt{det(g_{ij}(x,v))}}{\sigma_{\mu}(x)} \right).
 \end{equation}
 \end{Def}
\begin{Def} 
The  rate of changes of the distortion along a geodesic $\eta(t)$ is called $S_{\mu}$\emph{-curvature} (or simply S-curvature): 
\begin{equation}\label{S-curvature def} 
  S_{\mu}(x,v):=\frac{d\tau_{\mu}(\eta(t),\dot{\eta}(t))}{dt}|_{t=0} , 
\end{equation}
   where $\eta(t)$ is the geodesic starting from $x$ with initial velocity $v$. \end{Def}
 Thereby,  $S_{\mu}(t)$ denotes $S_{\mu}(\eta(t),\dot{\eta}(t))= \tau'_{\mu}(t)$. The $S_{\mu}$-curvature can be expressed in local coordinates as follows:\vspace{-0.2cm}
   \begin{equation}
   S_{\mu}(x,y)= \dot{\partial}_i\, G^{i}(x,y) - y^{i}\, \partial_{i}\left(\log(\sigma_{\mu}(x))\right).
   \end{equation}
\begin{Pro}\label{tau vanishies iff F is Riem}\cite[Proposition 4.2]{Shenbook16}  Given a Finsler $\mu$-space $(M,F, \mu)$, the following three statements are equivalent: \emph{(i)} $F$ is Riemannian; \emph{(ii) }The {distortion} $\tau_{\mu}$ vanishes identically; \textit{(iii)}  The {distortion} $\tau_{\mu}(x,v)$ is function of $x$ only.  
\end{Pro}
\begin{flushleft}
\textbf{Gradient, Hessian and Laplacian in Finsler geometry}.
\end{flushleft}

Now let us recall the definitions of gradient, Hessian and Laplacian in Finsler setting and some relations between them.  For further details,  see \cite{A.A. Tamim,  universalcomp}.

It is known that,  if $F$ is a Finsler structure on $M$, then $F$ induces at each point $x\in M$  a Minkowski norm on $T_x M$. Also, $F^*$, the dual structure of $F$,  induces a Minkowski norm on $T_x^*M$.  That is, $F^*: T^*M \rightarrow \R^+$  is defined, for all $(x,\alpha) \in T^*M$, by $$F^*(x,\alpha) :=\sup\{ \alpha (\xi) \,:\, \xi \in I_xM \}. $$  The dual metric associated to $F^*$ is given by \vspace{-0.25cm}
  $$g^*_{ij}(x,\alpha) :=\frac{1}{2}\,\frac{\partial^2}{\partial \alpha^i \partial \alpha^j} \big( F^{*2}(x,\alpha) \big).$$   
The \textit{Legendre transformation} $J: TM \rightarrow T^*M$ associated with the Finsler structure $F$ is defined, for any point $x \in M$, by 
$$J(x,y)= g_{ij}(x,y)\, y^{i}\, dx^{j},\,\, \,\,\forall y \in T_{x}M_{0} \text{ and  }J(0)=0 .$$ Let  $J^*: T^*M \rightarrow TM$ defined by \vspace{-0.3 cm}
$$J^{*}(x,\alpha)= g^*_{ij}\big( x,\alpha \big)\, \alpha_{i}\, \frac{\partial}{\partial x^j},\,\, \forall \alpha \in T^{*}_{x}M_{0} \text{ and  }J^{*}(0)=0 ,$$
where  $g^*_{ij}(x,\alpha) := g^{ij} (J^*(\alpha))$.

\begin{Def}\cite[\S 3.2]{Shlec} The gradient of a differentiable function $f:M \rightarrow \R$  at a point $x \in M$, where $df(x) \neq 0$,   is defined by 
\vspace{-0.3 cm}\begin{equation}\label{graddef}
\nabla f(x)=J^*\big( x,df(x) \big)= g^*_{ij}\big( x,df(x) \big)\,\, \frac{\partial f(x)}{\partial x^i}\, \frac{\partial}{\partial x^j}.
\end{equation}
\end{Def}
\begin{flushleft}
$df(x)$ can then be written in the form
\end{flushleft}
\begin{equation}\label{grad}
df(x,v)= g_{\nabla f(x)}(\nabla f(x), v), \,\, \forall v \in T_{x}M.
\end{equation}

\begin{Rem}
Unlike the Riemannian gradient, the gradient $\nabla f(x)$ is nonlinear.  It should be noted that, when $df(x) = 0$, the gradient $\nabla f(x)$ is  defined to be zero.
\end{Rem}
\begin{Def}\cite{Shlec}
A smooth function $f: M \longrightarrow \mathbb{R}$ is called a Finsler distance if $F(\nabla f)=1$.
\end{Def}
A distance function $r$ defined on an open subset $\Omega$ of $(M,F)$ has some interesting geometric properties. Indeed, $\nabla r$ is a unit vector field on $\Omega$ and it induces a smooth Riemannian metric on $\Omega$ defined by \vspace{-0.4 cm} $$\hat{F}(x,v):= \sqrt{g_{\nabla r}(v,v)},\,\, \forall v \in TM.$$ Further more, $\hat{F}(\hat{\nabla} r)= F(\nabla r) =1$ by \cite[Lemma 3.2.2]{Shlec}.
\begin{Def}\cite{Shenbook16}
The Hessian $H(f)$ of a $C^2$ function $f$ is defined on the set $U_{f}:=\{x \in M \,|\, df(x)\neq 0\}$ by
\begin{equation}\label{Hess}
H(f)(X,Y):= XY(f)- \overline{\nabla}^{\nabla f}_{X} Y(f)= g_{\nabla f}(\overline{\nabla}^{\nabla f}_{X} \nabla f, Y),\quad \forall \, X,Y \in TU_{f} \subset TM,
\end{equation}
where $\overline{\nabla}$ is the Chern (Rund) connection.
\end{Def}
Here is another definition of the Hessian.
\begin{Def}\cite[\S 14.1]{Shlec} The
Hessian of a $C^2$ function $f$ is a mapping $D^2 f: TM \longrightarrow \mathbb{R}$ defined by\vspace{-0.5 cm} \begin{equation} D^{2}f(v):=\frac{d^2}{dt^2}\left(f \circ \eta \right)|_{t=0},\end{equation} where $\eta$ is an arbitrary geodesic with initial velocity $v$. \end{Def}

\begin{Rem}
These above two definitions of Hessian in Finsler geometry are not equivalent. However, they are equivalent in the Riemannian case. Moreover, the Hessian of a distance function $r$ defined on an open subset $\Omega$ of $(M,F)$ satisfies
 \begin{equation}\label{hess and tcurvature}
 D^2 r(v)= \hat{D}^2 r(v)-T_{\nabla r}(v), \,\, \forall v \in T_{x}\Omega,
 \end{equation}
where $\hat{D}^2 r$ is the Hessian of $r$ with respect to $\hat{F}$ and  $T_{\nabla r}(v)$ is the $T$-curvature of $(M,F)$ \cite[Lemma 14.1.1]{Shlec}. 
\end{Rem}
Unlike the Laplace-Beltrami operator in the Riemannian case, there are several notions of Laplacian in Finsler geometry; each of them has different properties. We refer to \cite{differentFLaplacian} and \cite{A.A. Tamim} for further information. We choose the Shen's Laplacian \cite{Shlec} to work with. 
\begin{Def}\cite[\S 14.1]{Shlec}
Let $(M,F,d \mu)$ be a Finsler $\mu$-space. 
 For a  $ C^2 $ function $ f$, the \emph{Shen's Laplacian} $\Delta f$ of $f$ is defined by $\Delta f=div_{\mu}(\nabla f)$, that is,
\begin{eqnarray}\label{shenlapdef}
\Delta f &=&\frac{1}{\sigma_{\mu}(x)} \,\partial_{k} \left[ \sigma_{\mu}(x) \; g^{kl}(x,\nabla f(x)) \; \partial_{l} f \right], \end{eqnarray} 
 where $\sigma_{\mu}(x)$ is the volume density of the volume form $d \mu$.
\end{Def}
\begin{Rem}
Shen's Laplacian is fully nonlinear elliptic differential operator of the second order which depends on the measure $\mu$ and it is defined on $U_{f}:=\{x \in M\, |\, df(x)\neq 0\}$ by $(\ref{shenlapdef})$, and to be zero on $\{x \in M \,| \,df(x)= 0\}$.
\end{Rem}
\textbf{Finsler Mean Curvature.}  Z. Shen has defined in \cite{curdisvol} the notion of mean curvature for hypersurfaces in $(M,F,\mu)$, where $d \mu = d V_{BH} $ is the Busemann-Hausdorff volume form. However, this definition can be used for an arbitrary  volume measure.  

Let $N $ be hypersurface of $M$. Suppose that $r$ is a Finsler distance defined on an open subset $U$ of $M$, i.e., $F(\nabla r) =1$, such that $r^{-1}(s)= N \cap  M$ for some $s \in \text{Im}(r) \subset \mathbb{R }$. Let $d \nu _{t}$ be the induced volume form by $d \mu$ on $N_{t}:=r^{-1}(t)$. Let $c(t)$ be an integral curve of $\nabla r$ which starts from $c(0) \in N_{s}$. Thereby, for small $\epsilon >0$, the flow $\phi _{\epsilon}$ of $\nabla r$ satisfies $\phi _{\epsilon} (c(s)):= c(s + \epsilon)$. Thus, $\phi _{\epsilon}  :r^{-1}(s) \rightarrow  r^{-1}(s + \epsilon)$. Therefore, the pull-back $(n-1)$-form $(\phi _{\epsilon})^{*} d \nu _{s+ \epsilon}$ is a multiple of $ d \nu _{s}$ \cite[\S 14.3]{Shlec}. Hence, there exists a function $\Theta(x,\epsilon)$ on N such that 
$$(\phi _{\epsilon})^{*} d \nu _{s+ \epsilon} =\Theta(x,\epsilon)\, d \nu _{s}|_{x} ,\, \forall x  \in N.$$
It should be noted that, $\Theta(x,0):=1,\,\forall x  \in N.$ 
The Finsler mean curvature of the level hypersurface $N$ at $x$  with respect to $\nabla r(x)$ is defined as follows, \cite[\S 14.3]{Shlec},  
\begin{equation}
\Pi_{\nabla r}(x):= \frac{\partial}{\partial \epsilon}\,\log(\Theta(x, \epsilon))|_{\epsilon =0}  .
\end{equation}
In a special local coordinate system $(t, x^{a})$ in $M$ such that $\nabla r = \frac{\partial}{\partial t}$ and $d \mu = \sigma (t, x^a )\, d t \wedge d  x^a ,\, a=2,...,n-1 $, the function $\Theta (x, \epsilon)$ can be  expressed as
$$\Theta (x, \epsilon)= \frac{\sigma (s+ \epsilon, x^{a})}{\sigma (s, x^{a})},\,\,\forall x \in N. $$
\begin{Def}\cite[\S 14.3]{Shlec} 
The \emph{Finsler mean curvature} of the level hypersurface $r^{-1}(t)$ at $x \in M$  with respect to $\nabla r_{x}$ is defined by  
\begin{equation}\label{Fmeancurvaturedef}
\Pi_{\nabla r}(x):= \frac{d}{dt}\,\log(\sigma_{x}(t,x^{a}))|_{t=s} = \frac{\partial}{\partial \epsilon}\,\log(\Theta(x, \epsilon))|_{\epsilon =0}.
\end{equation}
\end{Def}
\begin{Lem}\cite[Proposition 14.3.1]{Shlec}
The Finsler Laplancian of a distance function $r$ satisfies $$\Delta \, r(x)=\Pi_{\nabla r}(x).$$
\end{Lem}
\begin{Lem}\cite[Lemma 5.1]{Shenbook16}
 The relation between the Hessian of a $C^2$ function $f$ defined on $U_{f}$ and it Laplacian is given by
\begin{equation}\label{tarce}
\Delta f  = tr_{g_{\nabla f}} H(f) - S(\nabla f).
\end{equation}
\end{Lem}
\begin{Rem}
Equation $(\ref{tarce})$ shows that, the Finsler Laplacian can not be viewed in general as the trace of Hessian. However, when the metric is Riemannian, it can be expressed as a trace of Hessian. Indeed, $g_{\nabla f}$ is the induced Riemannian metric on the open subset $U$ of $M$.  In other words, the following expression gives the relation between the mean curvature $\Pi_{\nabla r}(x)$  of the level hypersurface $r^{-1}(t)$ in $(U,F)$ and the mean curvature $\hat{\Pi}_{\nabla r}(x)$ of $r^{-1}(t)$ in $(U,g_{\nabla r})$ 
\begin{equation}\label{Finslermeancur}
\Pi_{\nabla r}(x)=  \hat{\Pi}_{\nabla r}(x) -S\left(\nabla r(x)\right),
\end{equation} which is equivalent to
\begin{equation}\label{FinslerLac.}
\Delta \, r(x)=   \hat{\Delta} \, r(x) -S\left(\nabla r(x)\right).
\end{equation}
 \end{Rem}
\section{Properties of Normal and Mean Curvatures of Geodesic Spheres}
In this section, we discuss  some properties of the Finsler normal and mean curvature of geodesic spheres of different radii. The Finsler mean curvature is the mean curvature  of the level hypersurfaces $r^{-1}(t)$ in $(U,F)$, where $U$ is open subset of M. 
Let $S_{p}(r):= \exp \left[ \textbf{S}_{p}(r)\right]= \{ x\in M\, | \,d_{F}(p,x)=r \}$ be the forward geodesic sphere cf.~\cite[Chapter 6]{BCS}. 
Here, we shall use the definition and properties of the shape operator found in \cite[\S 14.4]{Shlec}. For further reading, we refer to \cite{curdisvol, {Shenbook16}}.\\

Berwald manifolds have many characterizations. One of them is  that $(M,F)$ is Berwaldian if and only if its $T$-curvature vanishes \cite[Proposition 10.1.1]{Shlec}.
\begin{Rem}\cite[p.~128]{BCS}
Let $\sigma : [0,r] \longrightarrow (M,F)$ be a Finslerian geodesic, then its reverse $\gamma(s):=\sigma (r-s)$ is again a geodesic if one of the following conditions is satisfied:
\begin{itemize}
\item[(1)] The Finsler structure is of Berwald type;
\item[(2)] The Finsler structure is reversible. \end{itemize} 
\end{Rem}
The following result generalizes \cite[Proposition 2.1]{Shah03} from Riemannian to Berwald spaces.

\begin{Lem}\label{nor cur Ber}
Let $(M,F,d \mu)$ be a forward complete, simply connected  Berwald $\mu$-manifold without conjugate points. Let $\eta_{v}$ be the minimal unit speed geodesic such that $\eta_{v}(0)=x,\, \dot{\eta}_{v}(0)=v$. Then for all $t>0$, the family of the Finsler normal curvatures $\{\Lambda_{t}\}_{t}$ of the forward geodesic spheres $S_{t}(\eta_{v}(t))$ at $x$ with respect to the outward pointing normal vector is strictly decreasing with~$t$.
\end{Lem}
\begin{proof}
Let $\eta_{v}$ be the minimal unit speed geodesic such that $\eta_{v}(0)=x,\, \dot{\eta}_{v}(0)=v,\,\eta_{v}(r)=p,\,\eta_{v}(R)=q$, where $r < R$ are positive numbers. 
Consider two forward geodesic spheres $S_{r}(p),\,S_{R}(q)$ touching each other internally at $x$. The unit outward pointing normal vector field to $S_{r}(p)$ is $ \frac{\partial}{\partial r}:=\nabla r = -v$.
 The induced second fundamental form is $ \hat{L}_{r}: Span\{v^{\perp}\} \longrightarrow  Span\{v^{\perp}\}$ defined by $\hat{L}_{r}(y):= \overline{\nabla}_{y} v, $ where $\overline{\nabla}$ is the Chern connection. Assume that $J_{1},\,J_{2}$ are Jacobi fields along $\eta_{v}$ such that $J_{1}(0)=J_{2}(0)=y,\,J_{1}(r)=J_{2}(R)=0$ and $X$ is a piecewise  $C^{\infty}$ vector field along $\eta_{v}$ over  $[0,R]$, defined by 
\begin{align} 
\operatorname{X}(t)=
 \left\{ 
        \begin{array}{cc}
J_{1}(t), & \quad  \ 0\leq t \leq r; \\
0, & \quad  \ r\leq t \leq R .
        \end{array} \right.
\end{align}
 
 Now applying the \textit{Index Lemma} \cite[Lemma 7.3.2]{BCS} to $J_{2}$ and $X$, we get
 $$I(X,X) \leq I(J_{2},J_{2}), $$
 where $I(.\,,.)$ is the index form along $\eta$.
 Using the formula  \cite[Eqn.~(7.2.4)]{BCS} $$I(J_{2},J_{2})= g_{\dot{\eta}_{v}}(J_{2}',J_{2})|_{0}^{R}= - g_{\dot{\eta}_{v}}(J_{2}'(0),y).$$
 Similarly, $I(J_{1},J_{1})=- g_{\dot{\eta}_{v}}(J_{1}'(0),y) = I(X,X) .$ Therefore,
 \begin{equation}\label{R2funfrom}
 g_{\dot{\eta}_{v}}(J_{1}'(0),y) >  g_{\dot{\eta}_{v}}(J_{2}'(0),y).
 \end{equation}
 As the Chern connection is torsion free, we get $$J_{1}'(0):= \overline{\nabla}_{\dot{\eta}_v} J_{1}(0)= \overline{\nabla}_{\dot{\eta}_v} y= \overline{\nabla}_{y} \,{\dot{\eta}_v (0)} = \hat{L}_{r}(y).$$ Hence, it follows from (\ref{R2funfrom}) that
 \begin{equation}\label{Rnorcur}
 g_{\dot{\eta}_{v}}(\hat{L}_{r}(y),y) >  g_{\dot{\eta}_{v}}(\hat{L}_{R}(y),y),\,\,\, \forall y \in span\{v^{\perp}\},\, r<R.
 \end{equation}
Using \cite[Lemma 14.4.1]{Shlec}, we deduce that (\ref{Rnorcur}) is equivalent to saying that the induced normal curvature satisfies 
\begin{equation}\label{Riemnormalcur}
\hat{\Lambda}_{\nabla r}(y) >\hat{\Lambda}_{\nabla R}(y),\,\,\, \forall y \in span\{v^{\perp}\},\, r<R.
\end{equation}
According to  Eqn.~$(\ref{hess and tcurvature})$ along with the vanishing of the $T$-curvature in the Berwald case and taking into account the fact that the normal curvature $\Lambda_{t}(y)$ is equal to the Hessian,  we conclude that $D^2_{t}(y)= \hat{D}^2_{t}(y)$, which is equivalent to $\Lambda_{t}(y)= \hat{\Lambda}_{t}(y)$.
\end{proof}
\begin{Cor}\label{mean cur Berwald}
Under the assumptions of Lemma\emph{ \ref{nor cur Ber}}, the mean curvature of the forward spheres $S_{t}(\eta_{v}(t))$ is strictly decreasing with $t$.
\end{Cor}
\begin{proof}
Taking the trace of (\ref{Riemnormalcur})  gives
\begin{equation}\label{Riemmeancur}
{\Pi}_{\nabla r}(y) >{\Pi}_{\nabla R}(y),\,\,\, \forall y \in span\{v^{\perp}\},\, r<R.
\end{equation}
Hence, the result follows.
\vspace{-5.5 pt}
\end{proof}
In the view of \cite[Proposition 10.1.1]{Shlec}, the Finsler spaces with non-vanishing $T$-curvature are non-Berwalian. We now obtain a more general result  using the same technique  of the proof of Lemma \ref{nor cur Ber}.
\begin{Pro}
Let $(M,F,d \mu)$ be a forward complete Finsler $\mu$-manifold.  Let $\eta_{v}$ be the minimal unit speed geodesic such that $\eta_{v}(0)=x,\, \dot{\eta}_{v}(0)=v$.
If the $T$-curvature $T_{t}(y)$ is an increasing function in $t$, then for all $t>0$, the family of the Finsler normal curvatures $\Lambda_{t}$ of the forward geodesic spheres $S_{t}(\eta_{v}(t))$ at $x$ with respect to the outward pointing normal vector is strictly decreasing with t.
\end{Pro}
\begin{proof}
Using the same technique of proof Lemma \ref{nor cur Ber}, we get $\hat{\Lambda}_{t}(y)$ is strictly decreasing with $t$. Hence, Eqn.~$(\ref{hess and tcurvature})$, when $T_{t}(y)$ is a function increasing in $t$, implies that $\hat{\Lambda}_{t}(y)-T_{t}(y)$ is decreasing. Therefore, ${\Lambda}_{t}(y)$ is decreasing in~$t$.
\end{proof}

\begin{Rem}
One can easily see that Corollary \ref {mean cur Berwald} holds for forward complete Finsler harmonic $\mu$-manifold $(M,F,d \mu)$  with vanishing $S$-curvature. It is very useful to apply  Corollary \ref {mean cur Berwald} for some non-Berwald metrics with vanishing $S$-curvature. Examples of such metrics are:
\begin{itemize}
\item[(a)] Shen's fish tank metric \cite{Shenfish}.
\item[(b)] Non-Berwaldian Randers metrics \cite{ShenGeometricmeaning} with vanishing mean Berwald curvature $E$.
\item[(c)] Einstein Kropina metrics with respect to the Busemann-Hausdorff volume form \cite[Remark p.~313]{Shenbook16}.
\end{itemize} 
\end{Rem}
In the following we shall discuss the sign of the Finsler mean curvature. As it will be shown, this sign depends on the vanishing of the $S$-curvature.
\begin{Lem}\label{postivemeancurv}
A forward complete simply connected $(M, F, \mu)$ Finsler $\mu$-manifold without conjugate points with vanishing $S_{\mu}$-curvature has non-negative mean curvature. 
\end{Lem}
\begin{proof}
 When $S_{\mu}$-curvature vanishes, the Laplacian of a distance function at a point  $x \in M$ is given by $$\Delta r (x)= \Pi_{\nabla r}(x)  = \tr_{g_{\nabla r}} H(r) (x),$$ which follows from $(\ref{tarce})$. Let $ \eta:  [ 0, r(p) ] \longrightarrow M $ be a normal minimal geodesic joining $p$ and $x$. Therefore, $\dot{\eta} (r(x))= \nabla  r (x)$. Assume that $J_{1}, ...,J_{n-1}$ are the normal Jacobi fields along $\eta$ with $J_{i}(0)= 0$ and $J_{i}(r(x))= e_i$, where $\{\nabla  r, e_{i}\}^{n}_{i=1}$ is an orthonormal basis for $T_{x}M$ with respect to $g_{\nabla  r}$.
 For $x \in M$, using $(\ref{Hess})$, we have \vspace{-0.3 cm}
 \begin{eqnarray*}
 \tr_{g_{\nabla r}} H(r) (x) &=&  g^{ij}_{\nabla r}H(r) (e_{i}, e_{j}) = \sum^{n-1}_{i=1} H(r) (J_{i}, J_{i})|_{x} \\
 &=& \sum^{n-1}_{i=1} g_{\nabla r}(\overline{\nabla}_{J_{i}}^{\nabla r} \nabla r, J_{i})= \sum^{n-1}_{i=1} g_{\nabla r}(\overline{\nabla}_{\nabla r}^{\nabla r} J_{i}, J_{i})\\
 &=&\sum^{n-1}_{i=1} I(J_{i}, J_{i}),
 \end{eqnarray*}
 where $\overline{\nabla}$ is the Chern connection, which is torsion free, thereby the orthogonal vectors $\nabla r, \, \{J_{i}\}^{n-1}_{i=1}$ satisfy $\overline{\nabla}_{J_{i}}^{\nabla r} \nabla r = {\overline{\nabla}}_{\nabla r}^{\nabla r} J_{i}$.
Since $M$ is without conjugate points, one can apply  \cite[Proposition 7.3.1]{BCS} to get $$ I(J_{i}, J_{i}) \geq 0;\,\,\, \forall \, 1 \leq i \leq n-1. $$ Hence, $\tr_{g_{\nabla r}} D^{2}(r) \geq 0$.  Therefore, $\Pi_{\nabla r}(x)$ is non-negative.
\end{proof}
It should be noted that Lemma \ref {postivemeancurv} generalizes the corresponding result \cite[Proposition 2.2]{Shah03} from Riemannian to Finsler spaces with zero $S$-curvature.

\begin{Pro}
Let $(M, F, \mu_{BH})$  be a forward complete simply connected  Berwald $\mu_{BH}$-manifold without conjugate points. Then,
 $ \displaystyle\lim_{r \to \infty}  \Pi_{\nabla r}(x)$ exists and is non-negative. 
\end{Pro}
\begin{proof}
It follows from Corollary \ref{mean cur Berwald} and Lemma \ref{postivemeancurv}.
\end{proof}
\begin{Rem}\label{OhtaScurzero}
Ohta in \cite{OhtaSzero} showed that for an $n$-dimensional Randers metric $F = \alpha + \beta$, if there is a volume element $d\mu$ such that $S_{\mu}=0$, then $\beta$ is a Killing form whose  length with respect to $\alpha$ is constant. On the other hand, Z. Shen, in \cite{Shlec}, showed that the converse is also true. 
\end{Rem}
Consequently, we obtain the following result.
\begin{Pro}
Let  $(M, F:= \alpha + \beta,\, \mu_{BH})$ be a  forward complete simply connected Finsler $\mu_{BH}$-manifold of Randers type without conjugate points such that $\beta$ is a Killing form whose  length with respect to $\alpha$ is constant. Then, it has a non-negative mean curvature. 
\end{Pro}
\begin{proof}
It follows directly from Lemma \ref{postivemeancurv} and Remark \ref{OhtaScurzero}.
\end{proof}

\section{Harmonic Finsler manifolds}
In this section, we generalize  several kinds of harmonic manifolds to the Finsler setting. In order to do so, we need to recall polar coordinate system in Finsler geometry.\\

Actually, the Hopf-Rinow theorem shows that for a connected Finsler space, forward completeness is  equivalent to the fact that the exponential map $\exp_{x}$ is defined on the whole $T_{x}M$ \cite[Theorem 6.6.1]{BCS}. Thus, we will assume from now on that our Finsler manifold is forward complete.
We  follow \cite[\S 2.4.3, \S 7.1.1]{Shenbook16} in defining the Finsler polar coordinate system. \\

The polar coordinate system $(\textbf{r},\textbf{y})$ on each tangent space $T_{x}M_{0}$ with the Minkowskian norm $ F(x,.)$ is given, for all $u \in T_{x}M_{0}$,  by
$$\textbf{r}(u):=F(x,u),\,\,\textbf{y}:= \frac{u}{\textbf{r}(u)} \in I_{x}M.$$ 
Then the Finsler metric $g_{x}:=g_{ij}(x,y) \,dy^i \otimes dy^j$ at $x$ is given by
$$ g_{x} = d\textbf{r} \otimes d\textbf{r} + {\textbf{r}}^2 \dot{g}_{x};\,\,\, \dot{g}_{x}= \dot{g}_{ij} \,d \overline{\theta}^{i} \otimes \overline{\theta}^{j} \,\,\, i,j=1,...,n-1, $$ 
where $\{ \overline{\theta}^{j}\}_{j=1,...,n-1}$ is the spherical coordinates on $I_{x}M$ and $\dot{g}_{x}$ is the restriction of $g_{x}$ on $I_{x}M$.
 
 Let  $D_{x}:= M-Cut_{x}$, where $Cut_{x}$ is cut locus of $x$. It is clear that $U:=D_{x} - \{x\}$ is the maximal homeomorphic domain of $\exp_{x}$.  The polar coordinate system on $U$, denoted $(r,y)$,  is given for all
$ x_{o} \in D_{x}$, by
$$r(x_{o}):= \textbf{r} \circ \exp^{-1}_{x}(x_{o});\,\, y(x_{o})= \textbf{y}\circ \exp^{-1}_{x}(x_{o}).$$
In other words,
$$ \frac{\partial}{{\partial r}}|_{(r,y)} =(d \,\exp_{x})_{ry}(y) ;\,\, \frac{\partial}{{\partial \theta^{i}}}|_{(r,y)} = (d\, \exp_{x})_{ry}(r \,\frac{\partial}{\partial \overline{\theta}^{i}}),$$
where ${\theta}^{i}(x_{o}) = \overline{\theta}^{i} \circ \textbf{y} \circ  \exp^{-1}_{x}(x_{o}) =  \overline{\theta}^{i} \circ y(x_{o}). $\\

A volume form $d\mu$ on $M$ in the polar coordinate system  can be expressed as $$d \mu = \sigma_{x}(r,y) \,dr \wedge d \Theta,$$ where $d \Theta = d\theta^{1} \wedge ... \wedge d\theta^{n-1}$. Thus, it can be written in the form
 $d \mu =\overline{\sigma}_{x}(r,y) \,dr \wedge d \nu_{x}(y),$ where 
\begin{equation}\label{Fdensity}
\overline{\sigma}_{x}(r,y):= \frac{{\sigma}_{x}(r,y)}{\sqrt{det(\dot{g}(x,y))}},
\end{equation}
 and $d \nu_{x}(y)$ is the induced Riemannian volume form on   $(I_{x}M, \dot{g}(x,y))$ with respect to the induced Riemannian metric $\dot{g}(x,y)$. More detailed information about these coordinates can be found in~\cite{Shenbook16, volform}. 
 
\begin{Rem}\label{remak1}
$(i)$ One can see from the definition of $\overline{\sigma}_{x}(r,y)$ that it is like a compatibility condition that relates an arbitrary volume form with the Finsler structure. Furthermore, it generalizes the well known volume measures in Finsler geometry, namely  Busemann-Hausdorff and Holmes-Thomson volume forms. That is why, we choose definition (\ref{Fdensity}) to introduce harmonic manifolds in the Finsler framework. Even though there is no canonical measure in Finsler geometry like the volume measure in Riemannian geometry, as aforementioned, we will work with an arbitrary measure $\mu$ on M.

$(ii)$  It is known that, in contrast to the Riemannian case,  the volume of the indicatrix Vol$(I_x)$  varies as $x$ varies. However, for  Landsberg manifolds Vol$(I_x)$ is constant \cite[Theorem 2]{vlolandsberg}, i.e., the volumes of all unit tangent spheres  are equal to each other. 
 \end{Rem}
  
\begin{Def}
\label{remak1} Every function $f: \mathbb{R}^{+} \longrightarrow \mathbb{R}$ generates a \textit{radial function} $f_{x}$  around a point $x \in M$ defined by $f_{x}(x_{o})= f(r_{x}(x,x_{o}))$, where $r_{x}(x,x_{o})$ is the geodesic distance between $x, x_{o} \in M$ induced by the Finsler function.
\end{Def}
\begin{Def}\cite{BCS}
Let $\eta_{v} :[0,a] \rightarrow (M,F) $  be a unit speed geodesic emanating  
from a point $x \in M$ with initial velocity $v$. The forward injectivity radius  $i(x)$ at $x \in M$ is defined by  \[i(x)=\inf \{i_{x}(v)\, |\, v \in I_{x}M\}, \text{ where } i_{x}(v)= \sup \{s>0 \,|\, d_{F}(x,\eta_{v}(t))=t,\,\, \forall \, t \in (0,s)\}.\] The injectivity radius $i(M)$ of  $(M,F)$ is given by $i(M):= \inf \{i(x)\, |\, x \in M\}.$
\end{Def}

\begin{Rem}
It should be noted that, $f_{x}$ is well defined only for the points $x_{o}$ for which $r_{x}(x,x_{o})$ is less than the injectivity radius at x. When the injectivity radius of $M$ is infinity, $f_x$ is globally defined. 
\end{Rem}

In the following, we introduce and investigate  the notion of different types of harmonic manifolds in the Finsler setting. We  formulate these definitions for a forward complete Finsler manifold $(M, F, \mu)$ with an arbitrary measure $\mu$ on $M$. Similarly, one can define these notions for a backward complete or complete  $(M, F, \mu)$. 
\begin{Def}
A forward complete $(M, F, \mu)$ Finsler $\mu$-manifold  is called \emph{locally harmonic at p} $\in$ M if in polar coordinates the volume density function $\overline{\sigma}_{p}(r,y)$ is a radial function in a neighborhood of $p$. That is,  $\overline{\sigma}_{p}(r,y)$ is independent of $y \in I_{p}M$; thus it can be written briefly as  $\overline{\sigma}_{p}(r)$.  Moreover, when the injectivity radius of $M$ is infinity, $(M,F,\mu)$ is called \emph{globally harmonic} if in  polar coordinates the volume density function $\overline{\sigma}_{p}(r,y)$ is a radial function around each $ p \in M.$ 
\end{Def}
\begin{Lem}\label{sigma expression}
 The volume density function $(\ref{Fdensity})$ can be written in the form 
 \begin{equation}\label{voldenstydefFinsler}
\overline{\sigma}_{p}(r,y)=e^{-\tau(\dot{\gamma}_{y}(r))}\,\, det(A_{p}(r,y)):=e^{-\tau(\dot{\gamma}_{y}(r))}\, r^{n-1}\,\, \frac{\sqrt{det(g_{\frac{\partial}{\partial r}}(\frac{\partial}{\partial \theta^i}, \frac{\partial}{\partial \theta^j}))}}{\sqrt{det(\dot{g}(p,y))}},\end{equation}
 where ${\frac{\partial}{\partial r}}|_{(r,y)}= {(d\,\exp_{p})}_{ry}(y)=\dot{\gamma}_{y}(r)$ and $\dot{\gamma}_{y}(0)=y.$
\end{Lem}
\begin{proof}Let ${\gamma}_{y}(t):= \exp_{x}(ty)$ be the minimal geodesic in $(M,F)$ starting from $x$ in the direction of $y \in I_{x}M.$ 
The proof follows from Eqn.~$(\ref{distortion def})$ and applying the Gauss lemma.
\end{proof}
 The above result also appears in \cite[Lemma 3.1]{universalcomp}.
\begin{Rem}\label{modified har. def.}
Let $\psi_{p}(r,y):=e^{-\tau(\dot{\gamma}_{y}(r))} \,\, \frac{\sqrt{det(g_{\frac{\partial}{\partial r}}(\frac{\partial}{\partial \theta^i}, \frac{\partial}{\partial \theta^j}))}}{\sqrt{det(\dot{g}(p,y))}}$. Thus, the volume density can be expressed in the form
\begin{equation}\label{sigma=rphi}
\overline{\sigma}_{p}(r,y)=r^{n-1}\, \psi_{p}(r,y).
\end{equation}
Consequently, our definition of local and global harmonicity  reads: $(M,F,\mu)$ is harmonic if  $\psi_{p}(r,y)$ is a radial function.
\end{Rem}

\begin{Pro}
 Our definition of (local/global) harmonic Finsler manifold reduces to the exiting one in Riemannian geometry when the Finsler structure is Riemannian. 
\end{Pro}
\begin{proof}
It is clear, from Proportion \ref{tau vanishies iff F is Riem}, that the vanishing of $\tau$ is equivalent to the Finsler structure being  Riemannian.  Therefore, $\overline{\sigma}_{p}(r,y)=det(A_{p}(r,y))$ is independent of the chosen measure $\mu$. In other words,  $det(A_{p}(r,y))$  depends solely on the Riemannian metric $g$. In fact, it was proved in \cite[\S 4]{ShenAdv97} that, 
\begin{eqnarray*}\nonumber
det(A_{p}(r,y))&= &r^{n-1}\,\,det\left[g_{\frac{\partial}{\partial r}}\left(\left(d \,\exp_{p}\right)_{rv}\left({\frac{\partial}{\partial \theta^i}}\right),\left(d\, \exp_{p}\right)_{rv}\left({\frac{\partial}{\partial \theta^j}}\right) \right)\right].
\end{eqnarray*}
Hence,
\begin{eqnarray}\label{Riem density Jacobi field}
det(A_{p}(r,y))&= &r^{n-1} || J_{1} (t)\wedge ... \wedge J_{n-1}(t)||_{\frac{\partial}{\partial r}}, 
\end{eqnarray}
where $J_{1}, ...,J_{n-1}$ are Jacobi fields along $\gamma_{v}(t):= \exp_{x}(tv)$.
\end{proof}
The following result is a characterization of harmonic Finsler manifolds in terms of the mean curvature $\Pi_{\nabla r}(x)$.
\begin{Pro}\label{Charact}
A Finsler $\mu$-manifold $(M, F, \mu)$ is locally (globally) harmonic if and only if the Finsler mean curvature of all geodesic spheres of sufficiently small radii (all radii), expressed in polar coordinates, is a radial function.
\end{Pro}
\begin{proof}
Let $(M, F, \mu)$ be harmonic. This means that  $\overline{\sigma}_{x}(r,y)$ is a radial function and, consequently, the radial derivative of its logarithm is radial as well. Hence, $\Pi_{\nabla r}$ is a radial function.

 For the converse, let $x \in M$  and $\Pi_{x}(R,y)$ be the Finsler mean curvature of a forward geodesic sphere $S_{x}(R)$. Then, in the view of $(\ref{sigma=rphi}),$ we have
  $$\Pi_{x}(R,y)=\frac{d}{dr}\,\log(\overline{\sigma}_{x}(R,y))= \frac{n-1}{R} +  \frac{d}{dr}\,\log(\psi_{x}(R,y)).$$  
 
 Assume now that  $\Pi_{x}(R,y)$ is a radial function, i.e.  $\Pi_{x}(R,y)= \Pi (r(x,.))=\Pi_{x}(r)$, then
$$\Pi_{x}(r)- \frac{n-1}{r} = \frac{d}{dt}\,\log(\psi_{x}(t,y))|_{t=r}.$$
Solving this equation with the initial condition $\psi_{x}(0,y)=1$, yields
$$\log(\psi_{x}(r,y))- \log(\psi_{x}(0,y))=\int^{r}_{0} \left(\Pi_{x}(t) - \frac{n-1}{t} \right)dt .$$ Therefore,  $\psi _{x}(r,y)= e^{\int^{r}_{0} \left( \Pi_{x}(t) - \frac{n-1}{t} \right)dt }$ is a radial function. Hence, by Remark~\ref{modified har. def.} $(M, F, \mu)$ is harmonic.
\end{proof}

\begin{Cor}\label{Charact. Laplancain}
A Finsler $\mu$-manifold $(M, F, \mu)$ is  harmonic if and only if the Shen's Laplancian of a distance function  is a radial function. 
\end{Cor}
\begin{proof}
It follows directly from \cite[Proposition 14.3.1]{Shlec}, which states that the Shen's Laplancian of a distance $r$ satisfies $\Delta \, r(x)=\Pi_{\nabla r}(x)$.
\end{proof}

The flag curvature is a natural generalization of the sectional curvature. There are various characterisations for Finsler manifolds of constant flag curvature $K$. It is known that, the model Finsler spaces are not completely classified as in the Riemannian case. In general, there are infinitely many Finsler model spaces, which are not isometric to each other. For example, in the Finsler spaces of negative flag curvature $K$: the Funk metrics having $K=-\frac{1}{4}$  are forward complete and non-reversible Finsler metrics. However, the Hilbert metrics having $K=-1$ are complete and reversible Finsler metrics. One can find further information in \cite{BCS, {Shlec}, {Shenbook16}}. Generally, Finsler manifolds of constant flag curvature do not have constant $S_\mu$-curvature. For example, Bryant metrics on ${\mathbb{S}}^n$  have constant flag curvature $K=1$  and non-isotropic S-curvature, cf.~\cite{BCS, Shenbook16}.
\begin{Thm}\label{Thm1.1}
A forward complete Finsler $\mu$-manifold of constant flag curvature $K$  and constant $S_{\mu}$-curvature is globally harmonic.
\end{Thm}
\begin{proof}
Let $y \in I_{x}M$ be the initial velocity of a geodesic $\gamma_{y}(t)$. Assume that $J(t)$ be a Jacobi vector field along  $\gamma_{y}(t)$ with $J(0)=0$. As $(M,F)$ is a forward complete Finsler manifold of constant flag curvature $K$,  the Jacobi fields are  given by \cite[\S 9.7]{BCS}, $J(t)= \operatorname{\mathcal{S}}_{K}(t) \,E(t),$ where $E(t)$ is a parallel vector field along $\gamma_{y}(t).$ Therefore, by \cite[Eqn.~(3.5)]{universalcomp}, 
the volume density function \eqref{voldenstydefFinsler} has the following expression
\begin{equation}\label{constatflagden}
\overline{\sigma}_{p}(r,y)=e^{-\tau(\dot{\gamma}_{y}(r))}\,\mathcal{S}^{n-1}_{K}(r),
\end{equation} where
\begin{align}\label{solofjacobieq}
\operatorname{\mathcal{S}}^{n-1}_{K}(r) =
 \left\{ 
        \begin{array}{lll}
\frac{1}{\sqrt{K}}\,\sin^{n-1}\left(\sqrt{K}r \right) & \, \textrm{if} \ K > 0,\\
r^{n-1} & \, \textrm{if } \ K = 0, \\
\frac{1}{\sqrt{-K}}\,\sinh^{n-1}\left(\sqrt{-K}r\right) & \, \textrm{if} \ K < 0.
        \end{array} \right.
\end{align}
It is clear that when the $S_{\mu}$-curvature is constant, say $c \in \mathbb{R}$, the distortion is a radial function. Indeed, $\frac{d}{dt}\tau(\gamma_{y}(r),\dot{\gamma}_{y}(r))=c $ implies $\tau(r) = cr + c_{1}$, where $c_{1}$ is a constant. 
Hence, $\overline{\sigma}_{p}(r,y)$ is radial.
\end{proof}
\begin{Def}The Finsler mean curvature of horospheres  $\Pi_{\infty}$ is the Finsler mean curvature of spheres of infinite radius, which is defined by $\Pi_{\infty} =\displaystyle\lim_{r \to \infty}  \Pi_{\nabla r}(x)$. 
\end{Def}
\begin{Pro}
For a forward complete Finsler $\mu$-manifold of constant flag curvature $K$ and constant S-curvature $c$, the Finsler mean curvature of  forward geodesic spheres is a decreasing function in r. Furthermore, when $K \leq 0$ the Finsler mean curvature of horospheres is constant. 
\end{Pro}
\begin{proof}
Plugging (\ref{constatflagden}) in (\ref{Fmeancurvaturedef}), the Finsler mean curvature of a forward geodesic sphere in such spaces takes the from: 
\begin{align} \label{FmeancurcontKconstantS}
\operatorname \Pi_{\nabla r}(x) = -c +
 \left\{ 
        \begin{array}{lll}
(n-1)\sqrt{K}\,\cot\left(\sqrt{K}r \right) & \, \textrm{if} \ K > 0, \\
\frac{n-1}{r} & \, \textrm{if } \ K = 0, \\
(n-1)\sqrt{-K}\,\coth \left(\sqrt{-K}r \right) & \, \textrm{if} \ K < 0.
        \end{array} \right.
\end{align}
Taking the limit in (\ref{FmeancurcontKconstantS}) as $r \rightarrow \infty$, produces 
\begin{equation}\label{FmeancurhoroconstantS}
\Pi_{\infty}(x)=-c \quad (\text{for } K=0);\,\,\,\,\quad\Pi_{\infty}(x)=-c +(n-1) \sqrt{-K}\quad (\text{for } K<0).
\end{equation} 
\vspace*{-1.2cm}
\[\qedhere\]
\end{proof}
We would like to point out that $\Pi_{\infty}$  is an important geometric quantity in the study of asymptotic harmonic Finsler manifolds as will be shown later.
\begin{Cor}
Let $(M,F, \mu)$ be a $\mu$-Finsler manifold of  constant flag curvature. If the S-curvature is an increasing radial function, then the Finsler mean curvature of forward geodesic spheres is a decreasing function in $r$.
\end{Cor}
\begin{proof}
According to (\ref{constatflagden}), the Finsler mean curvature of  forward geodesic spheres in constant flag curvature spaces is given by
\begin{equation}\label{Finslermeancurconstant}
\Pi_{\nabla r}(x)=-\frac{d}{dr}\tau(\dot{\gamma}_{y}(r)) +\frac{d}{dr} \left(\log\left[ \mathcal{S}^{n-1}_{K}(r)\right] \right).
\end{equation}
 It is clear that  $\frac{d}{dr} \left(\log\left[ \mathcal{S}^{n-1}_{K}(r)\right] \right)$ is a decreasing function in $r$. Moreover, the S-curvature is an increasing radial function. Hence, $\Pi_{\nabla r}(x)$ is a decreasing function in $r$.
\end{proof}
\begin{Pro}
Let $(M,F, \mu_{BH})$ be a constant flag curvature Finsler manifold of a Randers type.   Then, $(M,F, \mu_{BH})$ is harmonic. 
\end{Pro}
\begin{proof}
 It is known that Randers metrics of constant flag curvature must have constant $S_{BH}$-curvature \cite{Shenbook16}. Then, the proof is completed by the use of Theorem \ref{Thm1.1}.
\end{proof}
Now, we provide some examples of our Theorem \ref{Thm1.1} for better understanding.
Assume that $(M,F,\mu)$ is a forward complete Finsler manifold with Busemann-Hausdorff volume measure. The following are examples of  globally harmonic Finsler manifolds which  have constant flag curvature $K$ and constant $S_{BH}$-curvature $c$.
\begin{itemize}
\item[\textbf{a.}] \textit{Minkowskian metrics}:
It is known that any Minkowskian metric has $K=0,\,S_{BH}= 0$. Therefore, \vspace{-0.5 cm}
\begin{equation}\label{Finsler mean cur of Minkoskian}
\Pi_{\nabla r}(x) = \frac{n-1}{r},\,\,\,\,\, \Pi_{\infty}(x)=0.
\end{equation}
This matches with the  examples of hypersurfaces in Minkowskian spaces $(\mathbb{R}^{n+1},F,\, \mu_{BH})$ given in \cite[\S 5]{Finsler2ndfform}. 

A notable example of a Minkowskian metric is the Berwald-Moor metric in $\mathbb{R}^n$ which is defined by $F(y)= (y^1 ... y^{n})^{\frac{1}{n}}.$
\item[\textbf{b.}]\textit{Shen's fish tank metric}: It is non-Berwald and non-projectively flat with $K=0,\,S_{BH}=0$ \cite{Shenfish}. Therefore, it is neither Riemmanian nor locally Minkowskian metric. The mean curvature of the geodesic spheres and the horospheres are given  respectively by \eqref{Finsler mean cur of Minkoskian}.

\item[\textbf{c.}]\textit{ Funk metrics}: They are projectively flat, \cite[Example 7.3.4]{Shlec}, and have $K=\frac{-1}{4},\,S_{BH}= \frac{n+1}{2}$. Thus, 
\begin{equation}\label{Finsler mean cur of Funk}
\Pi_{\nabla r}(x) = \frac{(n-1)}{2}\,\coth \left(\frac{r}{2}\right)- \frac{(n+1)}{2},\,\,\,\,\, \Pi_{\infty}(x)= -1.
\end{equation}
\item[\textbf{d.}] \textit{Bao and Shen constructed a family of non-Riemannian Finsler structures on odd-dimensional spheres}: the members of this family are non-projectively flat and have $K=1,\,S_{BH}=0$ \cite[Example 9.3.2]{Shlec}. Consequently, $$\Pi_{\nabla r}(x)=(n-1)\, \cot(r).$$
\end{itemize}
Unlike Riemannian harmonic manifolds \cite{Szabo}, we  prove the following.
\begin{Lem}
For a general harmonic Finsler manifold, the volume density function $\overline{\sigma}_{p}(r)$ depends on the starting point $p$. 
\end{Lem} 
\begin{proof}
This is due to the asymmetry of the Finsler distance.
\end{proof}  
However, for some special Finsler metrics, $\overline{\sigma}_{x}(r)$ is independent of $x$. 
\begin{Thm}\label{vdensityrenerseF}
For a globally harmonic reversible Finsler manifold, the volume density function $\overline{\sigma}_{p}(r)$ is independent of $p,$ for all $p \in M$.
\end{Thm}
\begin{proof}The Finsler structure is reversible if and only if the induced distance $d_{F}(p,q)$ is symmetric. We have  $\overline{\sigma}_{p}(r(q))= \overline{\sigma}_{p}(d_{F}(p,q))$ and
$\overline{\sigma}_{q}(r(p))= \overline{\sigma}_{q}(d_{F}(q,p))$. As, $$\overline{\sigma}_{p}(d_{F}(p,q)) =\overline{\sigma}_{q}(d_{F}(q,p)),\,\,\, \forall p,q \in M,$$ we conclude that $\overline{\sigma}_{q}(r)=\overline{\sigma}_{p}(r),\,\,\, \forall p,q \in M.$
\end{proof}
A geometric meaning of the zeros of the volume density function $\overline{\sigma}_{p}(r)$ in globally harmonic manifolds is given in the following result. 
\begin{Thm}\label{Blaschke}
For a globally harmonic Finsler manifold, the zeros of $\overline{\sigma}_{p}(r,y)$ are conjugate points of $p$.
\end{Thm} 
\begin{proof}
Let $\eta$ be a Finslerian geodesic joining $p,q \in M$. Assume that $\overline{\sigma}_{p}(r)=0$. That is,  $\overline{\sigma}_{p}(r_{p}(p,q))=0$.  Using Lemma \ref{sigma expression}, the vanishing of the volume density  $\overline{\sigma}_{p}(r)=0$ means that the exponential map is singular. Therefore, by \cite[Proposition 7.1.1]{BCS},  $p$ is conjugate to $q$. 
\end{proof}

\begin{Cor}\label{Pro. compact}
Let  $(M, F, \mu)$ be a locally harmonic Finsler manifold. If a conjugate point occurs, then $(M, F, \mu)$ is compact. 
\end{Cor}
\begin{proof}
Since $M$ is connected, then once a conjugate point $p$ occurs, it occurs everywhere on $S_{p}(r)$, which follows from $\overline{\sigma}_{p}(r,y)$ begin a radial function. Hence, every geodesic emanating from $p$ contains a cut point. Thus, $M$ is compact in the view of  \cite[Lemma 8.6.1]{BCS}. 
\end{proof}

\begin{Def}\cite{ Blaschke03, Blaschke19, ShiBanktesh}
A compact Finsler manifold $(M,F)$ is said to Blaschke if \[ i(M) = d(M),\] where  
  $$d(M):= \sup \{d_{F}(x,x')\, |\, x,\,x' \in M\} \text{ is the diameter of  } (M,F) .$$
\end{Def}
A consequence of Theorem \ref{Blaschke} is the next result:

\begin{Pro}
A compact locally harmonic Finsler manifold is a Blaschke Finsler manifold.
\end{Pro}
\begin{proof}
If a conjugate point occurs, then by Corollary \ref{Pro. compact} it occurs at the same distance with the same multiplicity  for every point $p \in M$. Hence, the proof is completed. 
\end{proof}
\begin{Pro}
Let $(M,F, \mu)$ be a complete simply connected Finsler manifold of constant flag curvature $K=1$. Then, $(M,F)$ is a Blaschke Finsler manifold.
\end{Pro}
\begin{proof}
As a  complete simply connected Finsler manifold of constant flag curvature $K=1$ is diffeomorphic to $S^n$ and all of its geodesics are closed with length of $2 \pi$  \cite{Shen1996}, then, it is a Blascke Finsler manifold.
\end{proof}
\begin{Pro}
Let $(M,g)$ be  a Blaschke Riemannian manifold. Assume that  $\beta$ is a closed $1$-form whose length $||\beta||_{g} <1$ and $\mu$ is either Busemann-Hausdorff, Holmes-Thompson or extreme volume measures on M.  
Then,  $(M, F:= g +\beta, \mu)$ is a Blaschke Finsler manifold.
\end{Pro}
\begin{proof}
 Given a Blaschke Riemannian manifold $(M,g)$, then all of whose geodesics are closed with the same length. Thus, the Randers metric $F:= g +\beta$ has reversible geodesics since the $1$-form $\beta$ is closed.  That is,
 all geodesics of $(M, F, d\mu)$ are closed with the same length as  those of $(M, g )$  \cite[Lemma 6.4]{Blaschke19}. Hence, it is a Blaschke Finsler manifold.
\end{proof}

Now, let us generalize another type of harmonic manifolds, namely, infinitesimal harmonic manifolds from the Riemannian to the Finsler setting.

\begin{Def}
$(M, F, \mu)$ is called an \emph{infinitesimal harmonic Finsler manifold} at $x \in M$ if it satisfies the condition:\\  $\forall n \in \mathbb{Z}^{+},\,\, \exists\, c_{n}(x) \in \mathbb{R}$ such that the radial derivatives of $\overline{\sigma}_{x}(r,y)$ at the origin is $c_{n}(x)$. \\
That is, $\forall n \in \mathbb{Z}^{+},\, \exists \, c_{n}(x) \in \mathbb{R}:$
\begin{equation}\label{IFatx}
 D^{(n)}_{Y_{x}}\overline{\sigma}_{x}(r,y)|_{r=0}= \frac{d^n}{dr^n}\Bigg( r \mapsto \overline{\sigma}_{x}(\exp_{x}(rY_{x}))\Bigg)(0) = c_{n}(x),\, \forall \, \, Y_{x} \in I_{x}M .
\end{equation}
\end{Def}
\begin{Def}
 $(M, F, \mu)$ is called an \emph{infinitesimal harmonic Finsler manifold} if it satisfies; $\forall n \in \mathbb{Z}^{+},\, \exists\, c_{n} \in \mathbb{R}$ such that 
$$ D^{(n)}_{Y}\overline{\sigma}_{x}(r,y)|_{r=0}  = c_{n},\, \forall \, Y \in IM.$$
\end{Def}
\begin{Rem}
The above definitions reduce to the  corresponding Riemannian ones when the Finsler metric is Riemannian \cite[Chapter 6, 6.26]{Besse}. Besse conjectured, in the Riemannian context, that infinitesimal harmonic at every point implies infinitesimal harmonic \cite[Chapter 6.C, 6.D]{Besse}. In the Finsler context, till the moment, we do not know the relation between infinitesimal harmonic at $x \in M$ and infinitesimal harmonic.
 \end{Rem}
 In \cite{curdisvol},  Shen proved the following Lemma \ref{Tay}  for Busemann-Hausdorff volume measure $\mu_{BH}$. However, we observe that it is true for any arbitrary volume measure $\mu$. This is because, the $S$-curvature  varies as the volume measure varies.
\begin{Lem}\label{Tay} Let $(M,F,\mu)$ be a forward complete Finsler $\mu$-space. The Taylor expansion of the volume density function of the forward geodesic sphere
$S_{x}(r)$ at $x\in M$ is given by 
\begin{equation}\label{Taylor}
\overline{\sigma}_{x}(r,y)= r^{n-1} \left \{1+S(y)\, r+\frac{1}{2} \left[-\frac{1}{3} \,Ric(y) + \dot{S}(y)+ S^{2}(y)\right] r^2 +...\right\}.
\end{equation} 
Therefore, the Finsler mean curvature of a geodesic sphere $S_{x}(r)$   at $c(t)$ is given by
\begin{equation}\label{Taylormean}
\Pi_{y}= \frac{n-1}{r} -S(y) -\frac{1}{3}\, Ric(y)\, r - \dot{S}(y) + O(r),
\end{equation}
where $S(\dot{c}(t))= S(y) + \dot{S}(y)+O(r),$ and
$c(t)$ is a geodesic with initial velocity $\dot{c}(0)= y \in I_{x}M$.
\end{Lem}
\begin{Rem}
Formula (\ref{Taylor}) shows that the Ricci curvature and S-curvature determine the local behavior of the measure of small metric balls around a point. In the Riemannian case, the coefficients of the Taylor expansion at zero of the volume density function are universal polynomials
in the curvature tensor and its covariant derivatives.
\end{Rem}
\begin{Thm}\label{constantEinsteinConstantScur}
Infinitesimal harmonic Finsler manifolds are  Einstein Finsler manifolds with constant S-curvature and constant Ricci curvature.
\end{Thm}  
\begin{proof}
The main idea is to use the Taylor expansion of volume density function at $x\in M$.
Applying the definition of infinitesimal harmonic Finsler manifold to (\ref{Taylor}), it follows that 
\begin{equation}S(y) = c_{1}, \,\,\,\left[-\frac{1}{3}\, Ric(y) + \dot{S}(y)+ S^{2}(y) \right] = c_{2}.\end{equation}Therefore, $S(y) = c_{1}, \,\,\,Ric(y)= -3\,\left[c_{2}-c_{1}^{2}\right]$. Hence the result. 
\end{proof}
\begin{Pro}\label{IF is GH}
An infinitesimal harmonic Finsler manifold is globally harmonic. The converse is true only when $M$ is analytic. 
\end{Pro}
\begin{proof}
In the view of Theorem \ref{constantEinsteinConstantScur} and  (\ref{Taylor}), we deduce that $\overline{\sigma}_{x}(r,y)$ must be a radial function. For the converse, assume that  
 $(M, F, \mu)$ be globally harmonic. This means that  its volume density function $\overline{\sigma}_{x}(r,y)$ is a radial function and, consequently, its derivative $ D^{(n)}_{Y_{x}}\overline{\sigma}_{x}(r)|_{r=0}$ is constant.  $M$ is  supposed to be analytic for the convergence of the Taylor expansion (\ref{Taylor}). 
\end{proof}
\begin{Thm}
Let $(M,F,\mu)$ be an infinitesimal harmonic Finsler manifold at every point $x$ in M. Then, $(M,F,\mu)$ is an Einstein Finsler manifold with isotopic S-curvature. 
\end{Thm} 
\begin{proof}
Applying (\ref{IFatx}) to $$\overline{\sigma}_{x}(r,Y_{x})= r^{n-1} \left\lbrace1+S(Y_{x})\, r+\frac{1}{2}\, \left[-\frac{1}{3}\, Ric(Y_{x}) + \dot{S}(Y_{x})+ S^{2}(Y_{x})\right] \,r^2 +...\right\rbrace,$$
we get, for all $x \in M$,
\begin{equation}S(Y_{x}) = c_{1}(x), \,\,\,\left[-\frac{1}{3}\, Ric(Y_{x}) + \dot{S}(Y_{x})+ S^{2}(Y_{x})\right] = c_{2}(x).\end{equation}Therefore, $S(Y_{x}) = c_{1}(x), \,\,\,Ric(Y_{x})= -3\,\left[c_{2}(x)-c_{1}^{2}(x) -\dot{c_{1}}(x)\right]$. Hence, both Ricci curvature and S-curvature are isotropic.
\end{proof}
It is known that Ledrappier \cite{Ledrappier} introduced asymptotic harmonic Riemannian manifolds which are considered as a generalization of harmonic Riemannian manifolds, in the sense that the mean value property of harmonic functions and Einstein condition are not known to hold \cite[Chapter 5]{harmonicbook16}.

Here, we  define \emph{asymptotic harmonic Finsler manifolds}, shortly \emph{AHF-manifolds}, and investigate some of their properties. 
\begin{Def}
A forward complete, simply connected Finsler $\mu$-manifold $(M, F, \mu)$ without conjugate points  is called an \emph{AHF-manifold} if the Finsler mean curvature of horospheres is a real constant $h$.
\end{Def}
 Consequently, a noncompact harmonic Finsler manifold with  constant Finsler mean curvature of horospheres is an AHF-manifold.
Examples of Finsler metrics where  $(M,F, \mu_{BH})$ being AHF-manifold are the following:
\begin{itemize}
\item[\textbf{a.}]
 Minkowskian metrics with $h= 0$, by \eqref{Finsler mean cur of Minkoskian}. 
\item[\textbf{b.}] Shen's fish tank metric $h= 0$, by \eqref{Finsler mean cur of Minkoskian}.
\item[\textbf{c.}] Funk metrics with the real constant $h= -1$, by  \eqref{Finsler mean cur of Funk}.
 \end{itemize}

Towards the investigation of AHF-manifolds, we prove the following results using \emph{Riccati equation} of the shape operator $\hat{L}_{r}$ induced by the Riemannian metric $g_{\nabla r}$ \cite[\S 14.4]{Shlec}.
\begin{Thm}\label{Thm6.7}
Let $(M,F, \mu)$ be an AHF-manifolds with constant S-curvature, say $c$. Then, $(M,F, \mu)$ has Ricci curvature bounded above by a constant that depends on the Finsler mean curvature of horospheres $h$ and $c$. 
\end{Thm}
\begin{proof}
The Riccati equation is given by
\begin{equation}
\frac{d}{dr}\,\hat{L}_{r}(v) + \hat{L}_{r}^{2}(v) + R_{v}=0.
\end{equation}
Thus, 
\begin{equation}\label{RiemRiccatieq}
\frac{d}{dt}\,\hat{\Pi}_{x_{t}} + \frac{1}{n-1}\,\hat{\Pi}_{x_{t}}^{2} +Ric(Y_{t})\, \leq 0.
\end{equation}
Substituting by $\Pi_{x_{t}}$ from (\ref{FinslerLac.}) into (\ref{RiemRiccatieq}), we get
\[\frac{d}{dt}\,\left\lbrace \Delta \, r(x_{t}) +S(\nabla r(x_{t})) \right\rbrace  + \frac{\left\{\Delta \, r(x_{t})+S(\nabla r(x_{t}))\right\}^2}{n-1}+ Ric(Y_{t}) \,\leq 0. \]
Since $\Delta \, r(x_{t})=h,\, S(\nabla r(x_{t}))=c$, the last inequality gives
\begin{equation}\label{Thm 6.7eq}
Ric(Y_{t}) \leq - \frac{(h +c)^2}{n-1}. 
\end{equation}
\vspace*{-1cm}
\[\qedhere\]
\end{proof}
In particular, for Berwald spaces with Busemann-Hausdorff measure, the S-curvature vanishes identically, therefore the bound is simpler. Precisely, 
\begin{equation}\label{RicboundBER}
Ric(Y_{t}) \leq -\frac{h^2 }{n-1} .
\end{equation} 

As a consequence of (\ref{RicboundBER}), we get information about $2$-dimensional AH-Berwald manifold. 
\begin{Cor}
An AH-Berwald manifold of dimension $2$ is either locally Minkowskian or Riemannian real hyperbolic space. 
\end{Cor}
\begin{proof}
Here, the inequality (\ref{RiemRiccatieq}) becomes
$\hat{L}_{r}^{2}(v) + R_{v}=0 \Longleftrightarrow h^2 + R_{v}=0.$ That is, $R_{v}=-h^2 $ which is constant. Now, applying Szabo's rigidity  result \cite[Theorem 10.6.2]{BCS}, one concludes that the Finsler structure $F$ is a locally Minkowskian metric when $h=0$ or a Riemannian metric when $h \neq 0$. In fact, the canonical Riemannian metric is a real hyperbolic metric.
\end{proof}
In the view of (\ref {Thm 6.7eq}), we get
\begin{Cor}
Let $(M,F, \mu)$ be an AHF-harmonic manifold satisfying
the hypothesis of Theorem \emph{\ref{Thm6.7}}. Then, Ricci curvature of $(M,F, \mu)$ is nonpositive.
\end{Cor}
\section{Harmonic Finsler manifolds of Randers type}\label{Harmonic Finsler manifolds of Rander type}
 One can consider a Randers metric as a modification of a Riemannian metric that leads to a particular Finsler metric. These metrics are an important class of special Finsler metrics for which many results were obtained, see for example \cite{BCS, Blaschke19, OhtaSzero}.  Furthermore, there is a way to find examples of harmonic Finsler manifolds which are of Randers type. 
\begin{Thm}\label{harmonicRanderthm}
Let $(M,\alpha)$ be a harmonic Riemmanian manifold. Let $\beta$ be a $1$-form such that its length $||\beta||_{\alpha}$ is a radial function and $||\beta||_{\alpha} <1$. Then, $(M,F:=\alpha+\beta, \mu)$ is  a harmonic Finsler manifold of Randers type, where $\mu$ is either a Busemann-Hausdorff, Holmes-Thompson or extreme volume measure on $M$.
\end{Thm}
\begin{proof}
As $(M,\alpha)$ is a harmonic Riemmanian manifold, the volume density function $\sqrt{det(\alpha_{ij})}$ of $dV_{\alpha}$ is a radial function, say $l(r)$. In other words, $dV_{\alpha}= l(r)\, dr \wedge d\Theta$.
Since, $||\beta||_{\alpha}$ is a radial function, then we have the following relations  \cite{volform}:
\begin{align}\label{FRanderdensty}\nonumber
dV_{HT}&=dV_{\alpha}=  l(r) dr \wedge d\Theta,\\ \nonumber
dV_{BH}&= (1- ||\beta||_{\alpha}^{2})^{\frac{n+1}{2}}\, dV_{\alpha},\\[-0.5cm] 
&&\\\nonumber
dV_{max} &= (1+||\beta||_{\alpha})^{n+1} dV_{\alpha},\\ \nonumber
dV_{min} &= (1-||\beta||_{\alpha})^{n+1} dV_{\alpha}.\nonumber
\end{align}
Hence, the corresponding volume density functions are radial functions. Consequently, $(M,F:=\alpha+\beta, \mu)$ is a harmonic  Finsler manifold of Randers type.
\end{proof}
\begin{Cor} We have the following sequence of inequalities of volume forms:
$$dV_{min} \leq dV_{BH} \leq dV_{HT} \leq dV_{max}.$$
\end{Cor}
\begin{proof}
This follows directly from (\ref{FRanderdensty}).
\end{proof}
\begin{Rem}\label{Jacobi fields in constant sec. curv.}
 In order to exemplify Theorem \ref {harmonicRanderthm}, we recall the well-known examples of harmonic Riemannian spaces. Using (\ref{Riem density Jacobi field}), one can calculate the  Riemannian volume densities of some known globally harmonic Riemannian manifolds,  cf.~{\cite[\S 6.A, 6.18]{Besse}}. Indeed, let ${\gamma}_{y}(t):= \exp_{x}(ty)$ be the minimal  geodesic in $M$ starting from $x$ in the direction of $y \in I_{x}M.$ For an orthonormal basis $\{e_i\}_{i=1}^{n-1}$ of ${y}^\perp$  and  $J_{1}, ...,J_{n-1}$ the normal Jacobi fields along $\gamma_y$ with $J_{i}(0)= 0$ and $\dot{J}_{i}(0)= e_i$. Then 
\begin{equation}\label{vol. den. table1}
 \psi_{p}(r,y)= || J_{1} (r)\wedge ... \wedge J_{n-1}(r)||_{\frac{\partial}{\partial r}}.
\end{equation} 

As we mentioned in the introduction, the Riemannian space forms, cf.~\cite{Pet16}, which are Riemannian manifolds of constant sectional curvature $\kappa$, are classified for only three canonical local  Riemannian metrics. First metric is on $\R^{n}$ when ($\kappa =0$), the second is on ${\mathbb{S}}^n $ when ($\kappa =1$) and the third is on $\mathbb{RH}^{n}$ when ($\kappa =-1$), up to scaling. Thus, the Jacobi fields are given by
\begin{align*} 
\operatorname{J}_{i}(t) =
 \left\{ 
        \begin{array}{lll}
\sin(t)\, E_{i}(t) & \quad \textrm{if} \ \kappa=1, \\
t\, E_{i}(t) & \quad \textrm{if } \ \kappa=0, \\
\sinh(t)\, E_{i}(t) & \quad \textrm{if} \ \kappa= -1,
        \end{array} \right.
\end{align*}
where $\{E_{i}(t)\}_{i=1}^{n-1}$  are the parallel extensions of $\{e_i\}_{i=1}^{n-1}$.
\end{Rem}
 The following table \cite{Shahth}, for normalized Riemannian metrics, can be written in view of Remark~\ref{modified har. def.} and (\ref{vol. den. table1}). 


\begin{center}
{Table 1: Riemannian volume densities}\\
\end{center}
\vspace{0.08in}
\noindent
\begin{center}
\begin{tabular}{|c|c|c|c|}
\hline
{Compact}&{Volume density}&{Noncompact}&
{Volume density}  \\
{harmonic}&{function}&{harmonic}&
{function}  \\ 
{manifold}&{}&{ manifold}&{}\\
\hline
      &&$\R^{n}$& $r^{n-1}$\\
\hline	  
${\mathbb{S}}^n$ & $\sin^{n-1}(r)$ & $\mathbb{R}\textbf{H}^{n}$& $\sinh^{n-1}(r)$ \\
\hline
$\mathbb{C}\textbf{P}^n$& $\sin^{2n -1}(r)\;\cos(r)$& $\mathbb{C}\textbf{H}^{n}$&
$sinh^{2n-1}(r)\;cosh(r)$\\
\hline
$\mathbb{H}\textbf{P}^n$& $\sin^{4n-1}(r)\;\cos^3 (r)$&
$\mathbb{H}\textbf{H}^n$& $\sinh^{4n -1}(r)\;\cosh^3 (r)$\\
\hline

$Ca\textbf{P}^2$ & $\sin^{15}(r)\;\cos^{7}(r)$&
$Ca\textbf{H}^2$&
$\sinh^{15}(r)\;\cosh^{7}(r)$  \\   
\hline
\end{tabular}
\end{center}
\noindent

\vspace{0.3 cm}
In the above table the symbol $\mathbb{C}\textbf{P}^n$ denotes complex projective space \cite{Szabo}, $\mathbb{H}\textbf{P}^n$ quaternionic projective space~\cite{Szabo},  $Ca\textbf{P}^2$ octonionic projective plane (Cayley projective plane) \cite{Szabo}, $\mathbb{R}\textbf{H}^{n}$ real hyperbolic space~\cite{RamRanj}, $\mathbb{C}\textbf{H}^{n}$  complex hyperbolic space~\cite{RamRanj},  $\mathbb{H}\textbf{H}^n$  quaternionic hyperbolic space~\cite{RamRanj},  $Ca\textbf{H}^2$ (or $\textbf{H}^2(\mathbb{O})$) complex Cayley hyperbolic plane \cite{Cay}.\\  For further information we refer to \cite{Besse, CharharRiem, {Shahth}}.   

\begin{Rem}
  It should be noted that, the first relation in \eqref{FRanderdensty} gives the Holmes-Thomson volume form $dV_{HT}$ of the Randers metric $F:=\alpha + \beta$ which coincides with the volume form $dV_{\alpha}$ of the Riemannian metric $\alpha$  (already written in Table 1). That is, Table 1 represents the Randers Holmes-Thomson volume densities which are the same as Riemannian volume densities. 
\end{Rem}
In the view of the above discussion,  the  following tables present examples of volume densities for some compact and noncompact harmonic Randers spaces $(M, F:= \alpha + \beta, d \mu)$, where $||\beta||_{\alpha} :=f(r) < 1$.
\newpage
\begin{center}
{Table 2: Randers Busemann-Hausdorff volume densities }\\
\end{center}
\vspace{0.1in}
\noindent
  \begin{center}
{\fontsize{10.5}{12}\selectfont{\begin{tabular}{|c|c|c|c|}
\hline
{Compact}&{Busemann-Hausdorff}&{Non}&
{Busemann-Hausdorff}  \\
{H.}&{volume density}&{-compact}&
{volume density}  \\ 
{spaces}&{function}&{ H. spaces}&{function}\\
\hline
      &&$\R^{n}$& $r^{n-1}\left[ 1-f^{2}(r)\right]^\frac{n+1}{2}$\\
\hline	  
${\mathbb{S}}^n$ & $\sin^{n-1}(r)\left[ 1-f^{2}(r)\right]^\frac{n+1}{2}$ & $\mathbb{R}\textbf{H}^{n}$& $\sinh^{n-1}(r)\left[ 1-f^{2}(r)\right]^\frac{n+1}{2}$ \\
\hline
$\mathbb{C}\textbf{P}^n$& $\sin^{2n -1}(r)\;\cos(r)\left[ 1-f^{2}(r)\right]^\frac{n+1}{2}$& $\mathbb{C}\textbf{H}^{n}$&
$sinh^{2n-1}(r)\;cosh(r)\left[ 1-f^{2}(r)\right]^\frac{n+1}{2}$\\
\hline
$\mathbb{C}\textbf{H}^n$& $\sin^{4n-1}(r)\;\cos^3 (r)\left[ 1-f^{2}(r)\right]^\frac{n+1}{2}$&
$\mathbb{H}\textbf{H}^n$& $\sinh^{4n -1}(r)\;\cosh^3 (r)\left[ 1-f^{2}(r)\right]^\frac{n+1}{2}$\\
\hline
$Ca\textbf{P}^2$ & $\sin^{15}(r)\;\cos^{7}(r)\left[ 1-f^{2}(r)\right]^\frac{n+1}{2}$&
$Ca\textbf{H}^2$&
$\sinh^{15}(r)\;\cosh^{7}(r)\left[ 1-f^{2}(r)\right]^\frac{n+1}{2}$  \\   
\hline
\end{tabular}}}
\end{center}	
\noindent
\vspace{0.3 cm}
\begin{center}
{Table 3: Randers maximum volume densities }\\
\end{center}
\vspace{0.1in}
\noindent
\begin{center}
{\fontsize{11}{11}\selectfont{\begin{tabular}{|c|c|c|c|}
\hline
{Compact}&{Maximum}&{Non}&
{Maximum}  \\
{H.}&{volume density}&{-compact}&
{volume density}  \\ 
{spaces}&{function}&{ H. spaces}&{function}\\
\hline
      &&$\R^{n}$& $r^{n-1}\left[ 1+f(r)\right]^{n+1}$\\
\hline	  
${\mathbb{S}}^n$ & $\sin^{n-1}(r)\left[ 1+f(r)\right]^{n+1}$ & $\mathbb{R}\textbf{H}^{n}$& $\sinh^{n-1}(r)\left[ 1+f(r)\right]^{n+1}$ \\
\hline
$\mathbb{C}\textbf{P}^n$& $\sin^{2n -1}(r)\;\cos(r)\left[ 1+f(r)\right]^{n+1}$& $\mathbb{C}\textbf{H}^{n}$&
$sinh^{2n-1}(r)\;cosh(r)\left[ 1+f(r)\right]^{n+1}$\\
\hline
$\mathbb{C}\textbf{H}^n$& $\sin^{4n-1}(r)\;\cos^3 (r)\left[ 1+f(r)\right]^{n+1}$&
$\mathbb{H}\textbf{H}^n$& $\sinh^{4n -1}(r)\;\cosh^3 (r)\left[ 1+f(r)\right]^{n+1}$\\
\hline
$Ca\textbf{P}^2$ & $\sin^{15}(r)\;\cos^{7}(r)\left[ 1+f(r)\right]^{n+1}$&
$Ca\textbf{H}^2$&
$\sinh^{15}(r)\;\cosh^{7}(r)\left[ 1+f(r)\right]^{n+1}$ \\   
\hline
\end{tabular}}}
\end{center}	
\noindent
\vspace{0.3 cm}

\begin{center}
{Table 4: Randers minimum volume densities }\\
\end{center}
\vspace{0.1in}
\noindent
\begin{center}

{\fontsize{11}{11}\selectfont{\begin{tabular}{|c|c|c|c|}
\hline
{Compact}&{Minimum }&{Non}&
{Minimum}  \\
{H.}&{volume density}&{-compact}&
{volume vensity}  \\ 
{spaces}&{function}&{ H. spaces}&{function}\\
\hline
      &&$\R^{n}$& $r^{n-1}\left[ 1-f(r)\right]^{n+1}$\\
\hline	  
${\mathbb{S}}^n$ & $\sin^{n-1}(r)\left[ 1-f(r)\right]^{n+1}$ & $\mathbb{R}\textbf{H}^{n}$& $\sinh^{n-1}(r)\left[ 1-f(r)\right]^{n+1}$ \\
\hline
$\mathbb{C}\textbf{P}^n$& $\sin^{2n -1}(r)\;\cos(r)\left[ 1-f(r)\right]^{n+1}$& $\mathbb{C}\textbf{H}^{n}$&
$sinh^{2n-1}(r)\;cosh(r)\left[ 1-f(r)\right]^{n+1}$\\
\hline
$\mathbb{C}\textbf{H}^n$& $\sin^{4n-1}(r)\;\cos^3 (r)\left[ 1-f(r)\right]^{n+1}$&
$\mathbb{H}\textbf{H}^n$& $\sinh^{4n -1}(r)\;\cosh^3 (r)\left[ 1-f(r)\right]^{n+1}$\\
\hline
$Ca\textbf{P}^2$ & $\sin^{15}(r)\;\cos^{7}(r)\left[ 1-f(r)\right]^{n+1}$&
$Ca\textbf{H}^2$&
$\sinh^{15}(r)\;\cosh^{7}(r)\left[ 1-f(r)\right]^{n+1}$ \\   
\hline
\end{tabular}}}
\end{center}	
\noindent

\vspace*{0.5 cm}

In Riemammian geometry, isoparametric hypersurfaces are a remarkable class of submanifolds studied by many geometers see, for example \cite[p.~87-96]{RIso}. On the other hand, the study of Finslerian isoparametric hypersurfaces has recently be started in \cite{isoparametric}. 
\begin{Def}\cite{isoparametric}
Let $(M,F,d \mu)$ be a forward complete Finsler $\mu$-space. A $C^{2}$ function $f: M\longrightarrow \mathbb{R}$ is called   \emph{isoparametric} in $(M,F,d \mu)$ if there is a smooth function
$a(t)$ and a continuous function $b(t)$ such that 
\begin{equation}\label{isoparametric fn. def.}
F(\nabla f)=a(f) ,\,\,\,\, \Delta f= b(f).
\end{equation}
Each regular level set $f^{-1}(t)$ is called
an \emph{isoparametric hypersurface} in M. 
\end{Def}

\begin{Pro}
A forward complete Finsler $\mu$-space is  harmonic if and only if the distance function $d_{F}$, induced by $F$, is isoparametric.
\end{Pro}
\begin{proof}
It is clear that a Finsler distance $d_{F}$ is a transnormal function as $F(\nabla d_{F})=1$. Corollary~\ref{Charact. Laplancain} states that  $(M,F,d \mu)$ is  harmonic if and only if the Laplacian of the distance function is a radial function. That is, $d_{F}$ satisfies $\Delta d_{F}= b(d_{F})$, for some function $b$. Thus, $d_{F}$ satisfies~\eqref{isoparametric fn. def.}. This complete the proof.
\end{proof}




\end{document}